\documentclass[10pt]{article}
\usepackage{amsmath}
\usepackage{fullpage}
\author{Thomas J. Haines,  Alexandra Pettet}
\title{Formulae relating the Bernstein and Iwahori-Matsumoto presentations
of an affine Hecke algebra}
\usepackage{amsfonts}
\usepackage{amssymb}
\newtheorem{theorem}{Theorem} [section]
\newtheorem{lemma} [theorem] {Lemma}
\newtheorem{proposition} [theorem] {Proposition}

\newtheorem{corollary} [theorem] {Corollary}

\newenvironment{proof} { \emph{Proof.} } { \rule{2mm}{2mm} \\}
\newenvironment{proof2} { \emph{Proof of Theorem \ref{SameTranslation}}. } {
\rule{2mm}{2mm} \\}
\newenvironment{proof3} { \emph{Proof of Theorem \ref{Drinfeld}}. } {
\rule{2mm}{2mm} \\}
\date{}
\begin{document}
\maketitle

\begin{abstract}
We consider the ``antidominant'' variants $\Theta^-_\lambda$ of the elements 
$\Theta_\lambda$ occurring in the Bernstein presentation of an affine Hecke algebra ${\mathcal H}$.  We find explicit formulae for $\Theta^-_\lambda$ in terms of the Iwahori-Matsumoto generators $T_w$ ($w$ ranging over the extended affine Weyl group of the root system $R$), in the case (i) $R$ is arbitrary 
and $\lambda$ is a {\em minuscule} coweight, or (ii) $R$ is attached to 
${\rm GL}_n$ and $\lambda = me_k$, where $e_k$ is a standard basis vector 
and $m \geq 1$.

In the above cases, certain {\em minimal expressions} for $\Theta^-_\lambda$ play a crucial role.  Such minimal expressions exist in fact for any coweight $\lambda$ for ${\rm GL}_n$.  We give a sheaf-theoretic interpretation of the existence of a minimal expression for 
$\Theta^-_\lambda$: the corresponding perverse sheaf on the affine Schubert 
variety $X(t_\lambda)$  is the push-forward of an explicit perverse sheaf on 
the Demazure resolution $m: \tilde{X}(t_\lambda) \rightarrow X(t_\lambda)$.  
This approach yields, for a minuscule coweight $\lambda$ of any $R$, or for an 
{\em arbitrary} coweight $\lambda$ of ${\rm GL}_n$, a conceptual albeit less explicit expression for the coefficient $\Theta^-_\lambda(w)$ of the basis element $T_w$, in terms of the cohomology of a fiber of the Demazure resolution.  

AMS Subject classification: 20C08, 14M15.

\end{abstract}

\section{Introduction}
Let $\mathcal{H}$ be the affine Hecke algebra associated to a root system.  
There are two well-known presentations of this algebra by generators and 
relations, the first discovered by Iwahori-Matsumoto \cite{I-M} and the 
second by
Bernstein \cite{Lusztig}, \cite{Lus2}; cf. 2.2.1 below.
The Iwahori-Matsumoto presentation
reflects the structure of the Iwahori-Hecke algebra $C^\infty_c(I\backslash 
G/I)$ of the split $p$-adic group $G$ attached to the root system: the 
generators $T_w$ correspond to the characteristic functions of Iwahori 
double cosets
$IwI$, where $w$ ranges over the extended affine Weyl group.   The Bernstein 
presentation reflects the description of the Hecke algebra as
an equivariant $K$-theory of the associated Steinberg variety, which plays a 
role in the classification of the representations of $\mathcal{H}$,
see \cite{KL}, \cite{CG}.  The Bernstein presentation has the advantage that 
one can construct a basis for the center of $\mathcal{H}$ by summing the 
generators $\Theta_\lambda$ over Weyl-orbits of coweights $\lambda$; the resulting functions are known as {\em Bernstein functions}.

It is of interest to give an explicit relation between the generators in 
these two
presentations.  More precisely, one would like to write each 
$\Theta_\lambda$ as an
explicit linear combination of the Iwahori-Matsumoto basis elements
$T_w$.  A direct consequence would be the explicit description of the
Bernstein functions (and thus the center of $\mathcal{H}$) in terms of the
Iwahori-Matsumoto basis.  This problem
was considered earlier by the first author \cite{Bernstein},\cite{TFSV} 
because of 
certain applications to the study of Shimura varieties, and was completely 
answered there for the case where $\lambda$ is a minuscule coweight.  More 
recently, O. Schiffmann \cite{Schiffmann} has given explicit formulae 
for all elements in a certain basis for the center $Z({\mathcal H})$ of an affine Hecke algebra
${\mathcal H}$ of type $A$; from this one can derive a formula for the Bernstein function $z_\mu$, where $\mu$ is any 
dominant coweight of a group of type $A$.

In this paper we consider the ``antidominant'' variants $\Theta^-_\lambda$ 
of the elements $\Theta_\lambda$.  The support of these functions is 
somewhat more regular than the original functions $\Theta_\lambda$, cf. 
Lemma 2.1.  In section 3 we consider the case where $\lambda$ is a minuscule 
coweight, and we prove the following explicit formula for 
$\Theta^-_\lambda$.

\begin{theorem}
Let $\lambda \in X_*$ be minuscule.  Then
$$
\Theta^-_\lambda = \sum_{\{ x \, : \, \lambda(x) = \lambda \}} 
\tilde{R}_{x,t_\lambda}(Q) \tilde{T}_x.
$$
\end{theorem}
Here $\lambda(x)$ is the translation part of $x$ ``on the left'' defined by the decomposition $x = t_{\lambda(x)}w$ $\,\,(w \in W_0)$, $\, \tilde{T}_x$ is 
a renormalization of the usual Iwahori-Matsumoto generator $T_x$, $Q = 
q^{-1/2} - q^{1/2}$, and $\tilde{R}_{x,y}(Q)$ is a variant of the usual 
$R$-polynomial of Kazhdan-Lusztig \cite{KL1}.

The formula above is analogous to the expression for $\Theta_\lambda$ found 
by the first author in 
Proposition 4.4 of \cite{TFSV}:
$$
\Theta_\lambda = \sum_{\{ x \, : \, t(x) = \lambda \}}
\tilde{R}_{x,t_{\lambda}}(Q) \tilde{T}_x.
$$
Here $t(x)$ is the translation part of $x$ ``on the right'' defined by the decomposition $x = wt_{t(x)}$ $\,\,
(w \in W_0)$.  However our proof
is simpler and more direct than that of loc.cit., and the same arguments 
appearing here also give a short proof of the formula for $\Theta_\lambda$.  
In fact one can derive the formula for $\Theta^-_{\lambda}$ from that for 
$\Theta_{-\lambda}$, and vice-versa.  Indeed, if 
$\iota: {\mathcal H} \rightarrow {\mathcal H}$ denotes the anti-involution determined by $q^{1/2} \mapsto q^{1/2}$ and 
$T_x \mapsto T_{x^{-1}}$, then 
$\iota(\Theta_{-\lambda}) = \Theta^-_{\lambda}$ and $\iota$ interchanges the formulae for $\Theta_{-\lambda}$ and $\Theta^-_{\lambda}$. 
We remark that 
Theorem 1.1 remains valid for Hecke algebras with arbitrary parameters.

In the fourth section we study coweights of ${\rm GL}_n$ of the form 
$\lambda = me_k$, where $e_k$ is the $k$-th standard basis vector and $m \in 
\mathbb{Z}_+$.  The case $m=1$, studied in \cite{Bernstein} and \cite{TFSV}, 
has relevance to a certain family of Shimura varieties with bad reduction, known as 
the {\em Drinfeld case}.
The general case is referred to as {\em multiples of the Drinfeld case}.  
We prove the following formula for 
$\Theta^-_{me_k}$.

\begin{theorem}
Let $1 \leq k \leq n$ and $m \geq 1$.  Then
$$
\Theta^-_{me_k} = \sum_{\{x \, : \, \lambda(x) \preceq me_k \}}
\tilde{R}_{x,t_{me_k}}(Q) \tilde{T}_x.
$$
\end{theorem}
Here $\preceq$ denotes the usual partial order on the lattice $X_*$.

Theorems 1.1 and 1.2 yield explicit expressions for the Bernstein functions 
$z_\mu$
($\mu$ minuscule) and $z_{me_1}$, respectively; see Corollary 3.6 and 4.2.  
The expressions in these special cases seem much simpler than the 
corresponding ones given by Schiffmann \cite{Schiffmann}.

Theorems 1.1 and 1.2 rely on the existence of certain 
{\em minimal expressions} for $\Theta^-_\lambda$: these are expressions of the form
$$
\Theta^-_\lambda = \tilde{T}^{\epsilon_1}_{t_1} \cdots \tilde{T}^{\epsilon_r}_{t_r}\tilde{T}_\tau,
$$
where $t_\lambda = t_1 \cdots t_r \tau$ ($t_i \in S_a, \, \tau \in \Omega$) is a reduced expression and $\epsilon_i \in \{ 1,-1 \}$ for every $1 \leq i \leq r$.  In the final two sections, we discuss how one can approach a general formula for $\Theta^-_\lambda$ when 
$\lambda$ is an arbitrary coweight of ${\rm GL}_n$, through minimal 
expressions (which always exist in this setting, cf. section 5).  The result 
is much less explicit than Theorems 1.1 and 1.2, and 
involves the geometry of the Demazure resolution $\widetilde{X}(t_\lambda) \rightarrow X(t_\lambda)$ of 
the affine Schubert
variety $X(t_\lambda)$.  We define a perverse sheaf 
$\Xi^-_\lambda$ on the affine flag variety whose corresponding function in the 
Hecke algebra is $\varepsilon_\lambda \Theta^-_\lambda$.  It turns out that $\Xi^-_\lambda$ is supported on 
$X(t_\lambda)$.  We see that the existence of a minimal expression for  
$\Theta^-_\lambda$  is analogous to the existence of a certain explicitly determined perverse sheaf on $\widetilde{X}(t_\lambda)$ whose push-forward to $X(t_\lambda)$ is $\Xi^-_\lambda$.  More precisely, we conclude the paper with the following result (cf. Theorem \ref{minimal_sheaf}, 
Corollary \ref{generalformula} for a completely precise statement).

\begin{theorem}
Let $\lambda$ be a minuscule coweight of a root system, or any coweight for ${\rm GL}_n$.  Choose a minimal expression for $\Theta^-_\lambda$ and let $m: \widetilde{X}(t_\lambda) \rightarrow X(t_\lambda)$ denote the corresponding Demazure resolution for $X(t_\lambda)$.  Then there exists an explicit perverse sheaf ${\mathcal D}$ on $\widetilde{X}(t_\lambda)$ (determined by the choice of minimal expression for $\Theta^-_\lambda$) such that 
$$
Rm_*({\mathcal D}) = \Xi^-_\lambda.
$$
Consequently, if we denote the coefficient of $T_x$ in the expression for $\Theta^-_\lambda$ by $\Theta^-_\lambda(x)$, then we have
$$
\Theta^-_\lambda(x) = \varepsilon_\lambda {\rm Tr}({\rm Fr}_q, H^\bullet(m^{-1}(x), {\mathcal D})),
$$
for any $x \leq t_\lambda$ in the Bruhat order.  Here the right hand side denotes the alternating trace of Frobenius on the \'{e}tale cohomology of the fiber over $x \in X(t_\lambda)$ with coefficients in the sheaf ${\mathcal D}$.
\end{theorem}
 
Given a coweight $\lambda$ for an arbitrary root system, let $\lambda_d$ denote the dominant coweight in its Weyl-orbit.  
We remark that there is a similar formula for $\Theta^-_\lambda$, {\em provided that 
$\lambda_d$ is a sum of minuscule dominant coweights}.

\section{Preliminaries}
\subsection{Affine Weyl group}
Let $(X^*,X_*,R,\check{R},\Pi)$ be a root system, where $\Pi$ is the set of
simple roots. The Weyl
group $W_0$ is generated by the set of simple reflections
$\{ s_\alpha : \alpha \in \Pi \}$.

We define a
partial order $\preceq$ on $X_*$ (resp. $X^*$) by setting
$\lambda \preceq \mu$ whenever $\mu - \lambda$ is a linear combination
with nonnegative integer coefficients of elements of
$\{ \check{\alpha} : \alpha \in \Pi \}$
(resp. $\{ \alpha : \alpha \in \Pi \}$).
We let $\Pi_m$ denote the set of roots $\beta \in R$ such that
$\beta$ is a minimal
element of $R \subset X^*$ with respect to $\preceq$.

In section 4 we will use the following description of the relation $\preceq$ 
for coweights of ${\rm GL}_n$: $(\lambda_1, \dots,\lambda_n) \preceq 
(\mu_1,\dots \mu_n)$ if and only if $\lambda_1 + \cdots + \lambda_i \leq 
\mu_1 + \cdots + \mu_i$ for $1 \leq i \leq n-1$, and $\lambda_1 + \cdots + 
\lambda_n = \mu_1 +
\cdots + \mu_n$.

Let $\widetilde{W}$ be the semidirect product
$X_* \rtimes W_0 = \{ t_x w : w \in W_0, x \in X_* \}$, with multiplication
given by $t_x w t_{x'} w' = t_{x+w(x')} ww'$.
For any $x \in \widetilde{W}$, there exists a unique expression
$t_{\lambda(x)} w$, where $w \in W_0$ and $\lambda(x) \in X_*$.

Let
$$
S_a =
\{ s_\alpha : \alpha \in \Pi \} \cup
\{ t_{-\check{\alpha}} s_\alpha : \alpha \in \Pi_m \} \subset \widetilde{W}.
$$
Define length
$l: \widetilde{W} \rightarrow \mathbb{Z}$ by
$$
l(t_x w) =
\sum_{\alpha \in R^+ : w^{-1}(\alpha) \in R^-}
|\langle \alpha,x \rangle -1| +
\sum_{\alpha \in R^+ : w^{-1}(\alpha) \in R^+}
|\langle \alpha,x \rangle|.
$$
Let $\check{Q}$
be the subgroup of $X_*$ generated by $\check{R}$. The subgroup
$W_a = \check{Q} \rtimes W_0$ of $\widetilde{W}$ is a Coxeter group with
$S_a$ the set of simple reflections. The subgroup is normal and admits a
complement
$\Omega = \{ w \in \widetilde{W} : l(w) = 0 \}$.

For $w \in \widetilde{W}$ denote $\varepsilon_w = (-1)^{l(w)}$ and $q_w = q^{l(w)}$ (for $q$ any parameter).

The Coxeter group $(W_a,S_a)$ comes equipped with the Bruhat order $\leq$.  
We extend it to $\widetilde{W}$ as follows: we say $w\tau \leq w'\tau'$
($w,w' \in W_a$, $\tau, \tau' \in \Omega$) if $w \leq w'$ and $\tau = 
\tau'$.

Let $\mu \in X_*$ be dominant. Following Kottwitz-Rapoport
\cite{KR}, we say $x \in \widetilde{W}$ is
$\mu$-\emph{admissible} if $x \leq t_{w(\mu)}$ for some $w \in W_0$.  We 
denote the set of $\mu$-admissible elements by ${\rm Adm}(\mu)$.

\subsection{Hecke algebra}
\subsubsection{Presentations}

The braid group of $\widetilde{W}$ is the group generated by $T_w$
($w \in \widetilde{W}$) with relations
$$
T_w T_{w'} = T_{ww'} \qquad
{\rm whenever} \ l(ww') = l(w) + l(w').
$$
The Hecke algebra $\mathcal{H}$ is defined to be the quotient of the group
algebra (over $\mathbb{Z}[q^{1/2},q^{-1/2}]$) of the braid group of
$\widetilde{W}$, by the two-sided ideal generated by the elements
$$
(T_s +1)(T_s -q),
$$
for $s \in S_a$. The image of $T_w$ in $\mathcal{H}$ is again denoted by
$T_w$. It is known that the elements $T_w$ ($w \in \widetilde{W}$) form a
$\mathbb{Z}[q^{1/2},q^{-1/2}]$-basis for $\mathcal{H}$.  The presentation of
$\mathcal{H}$ using the generators $T_w$ and the above relations is called 
the
{\em Iwahori-Matsumoto presentation}.

For any $T_w$, define a renormalization
$\tilde{T}_w = q^{-l(w)/2} T_w$.
Define an indeterminate $Q = q^{-1/2} - q^{1/2}$. The elements
$\tilde{T}_w$ form a basis for $\mathcal{H}$, and the usual relations can
be written as
\[\tilde{T}_{s} \tilde{T}_{w} =
\left \{
\begin{array}{ll}
\tilde{T}_{sw}, & {\rm if} \ l(sw) = l(w) +1, \\
-Q \tilde{T}_w + \tilde{T}_{sw}, & {\rm if} \ l(sw) = l(w) -1,
\end{array}
\right. \]
for $w \in \widetilde{W}$ and $s \in S_a$. There is also a right-handed
version of this relation.  Note that $\tilde{T}^{-1}_s = \tilde{T}_s + Q$.

We will denote $\tilde{T}_{t_\lambda}$ ($\lambda \in X_*$) simply by
$\tilde{T}_\lambda$.

For $\lambda \in X_*$, define
$$
\Theta_\lambda =
\tilde{T}_{\lambda_1} \tilde{T}^{-1}_{\lambda_2}
$$
where $\lambda = \lambda_1 - \lambda_2$, and $\lambda_1,\lambda_2$
are dominant.  The elements $\Theta_\lambda$ generate a commutative 
subalgebra of $\mathcal{H}$.  It is known that the elements $\Theta_\lambda 
T_w$
($\lambda \in X_*$, $w \in W_0$) form a $\mathbb{Z}[q^{1/2},q^{-1/2}]$-basis 
for $\mathcal{H}$.  These generators satisfy well-known relations (see Prop. 
3.6, \cite{Lus2}); in case the root system is simply connected, these are 
given by the formula
$$
\Theta_\lambda T_s - T_s\Theta_{s(\lambda)} =
(q-1)\frac{\Theta_\lambda - \Theta_{s(\lambda)}}
{1 - \Theta_{-\check{\alpha}}},
$$
where $s = s_\alpha$ and $\alpha \in \Pi$.  The presentation of 
$\mathcal{H}$ with generators $\Theta_\lambda T_w$ and the above relations 
is called the {\em Bernstein presentation}.

We also define
$$
\Theta^-_\lambda =
\tilde{T}_{\lambda'_1} \tilde{T}^{-1}_{\lambda'_2}
$$
where $\lambda = \lambda'_1 - \lambda'_2$, and $\lambda'_1,\lambda'_2$
are {\em antidominant}.

The involution $a \rightarrow \overline{a}$ of
$\mathbb{Z}[q^{1/2},q^{-1/2}]$ determined by
$q \mapsto q^{-1}$
extends to an involution
$h \rightarrow \overline{h}$, given by
$$
\overline{
\sum a_w T_w} =
\sum \bar{a}_w T^{-1}_{w^{-1}}.
$$

It is immediate that
$\overline{\Theta}_\lambda = \Theta^-_\lambda$.  Clearly the Bernstein 
presentation gives rise to an analogous presentation using the generators 
$\Theta^-_\lambda T_w$ in place of $\Theta_\lambda T_w$.

\subsubsection{Bernstein functions}

For each $W_0$-orbit $M$ in $X_*$, define the Bernstein function
$z_M$ attached to $M$ by
$$
z_M =
\sum_{\lambda \in M} \Theta_\lambda.
$$
When the $W_0$-orbit $M$ contains
the dominant element $\mu$, this function is denoted by $z_\mu$.

From Corollary~8.8 of Lusztig \cite{Lusztig}, we have
$z_\mu = \overline{z}_\mu$. Consequently,
$$
z_\mu =
\sum_{\lambda \in W_0(\mu)}
\Theta^-_\lambda.
$$

\subsubsection{A support property}

The preceding formula implies that when one studies Bernstein functions
there is no harm in working with the functions
$\Theta^-_\lambda$ instead of the functions $\Theta_\lambda$. We do so
in this paper because their supports enjoy a nice regularity property,
given by the following lemma.

\begin{lemma}\label{support}
For $\lambda \in X_*$, we have
$$
{\rm supp}(\Theta^-_{\lambda}) \subset
\{ x \ : \ \lambda(x) \preceq \lambda \}.
$$
\end{lemma}
\begin{proof}
Write
$$
\Theta^-_\lambda =
\sum_{y \leq t_\lambda}
a_y(Q) \tilde{T}_y,
$$
where $a_y(Q) \in \mathbb{Z}_+[Q]$
(see Lemma~5.1 and Corollary~5.7 of \cite{Bernstein}).

Choose a dominant coweight $\mu'$ such that $\mu' + \lambda(x)$ is
also dominant for any $x$ in the support of $\Theta^-_\lambda$.
Thus we have
$$
\tilde{T}^{-1}_{-(\mu' + \lambda)} = \Theta^-_{\mu' + \lambda} =
\Theta^-_{\mu'} \Theta^-_\lambda = \sum_y a_y(Q) \tilde{T}^{-1}_{-\mu'} 
\tilde{T}_y.$$
Let $y \in {\rm supp}(\Theta^-_\lambda)$.  We claim that $t_{\mu'}y$ 
belongs to 
the support of $\tilde{T}^{-1}_{-\mu'} \tilde{T}_y$.  Indeed, under the 
specialization map 
${\mathcal H} \rightarrow {\mathbb Z}[\widetilde{W}]$ determined by $q^{1/2} \mapsto 1$, the element 
$\tilde{T}^{-1}_{-\mu'} \tilde{T}_y$ maps to $t_{\mu'}y$.  
Since no cancellation occurs on the right hand side above, we see 
from this that
$t_{\mu'} y \in$ supp$(\tilde{T}^{-1}_{-(\mu' +\lambda)})$, and thus
$$
t_{\mu' + \lambda(y)} w_y =
t_{\mu'} y \leq
t_{\mu' + \lambda},
$$
where $y = t_{\lambda(y)}w_y$.  Since $\mu' + \lambda(y)$ and $\mu' + 
\lambda$ are both dominant, it is
well-known that this implies
$\mu' + \lambda(y) \preceq \mu' + \lambda$. The lemma follows.
\end{proof}

In the case where $\lambda$ is minuscule, this statement can be considerably 
sharpened; see Corollary \ref{support2}.

We remark that Lemma \ref{support} plays a key role in the proof of Theorem 
\ref{Drinfeld}.

\subsubsection{R-polynomials}

\noindent For any $y \in \widetilde{W}$, let
$y = s_1 \cdots s_r \tau$ $(s_i \in S_a, \tau \in \Omega$) be a reduced
expression for $y$.  Then for any $x$, we can write
$$
\tilde{T}^{-1}_{y^{-1}} =
\sum_{x \in \widetilde{W}}
\widetilde{R}_{x,y}(Q) \tilde{T}_x.
$$
where $\widetilde{R}_{x,y}(Q) \in \mathbb{Z}[q^{1/2},q^{-1/2}]$. These
coefficients $\widetilde{R}_{x,y}(Q)$ can be thought of as polynomial
expressions in $Q$ (as the notation suggests)
because of the identity
$$
\tilde{T}^{-1}_{y^{-1}} =
(\tilde{T}_{s_1} +Q) \cdots
(\tilde{T}_{s_r} +Q) \tilde{T}_{\tau}.
$$

\section{The minuscule case}
We say $\lambda \in X_*$ is {\em minuscule} if $\langle \alpha, \lambda 
\rangle \in \{0, \pm 1 \}$, for every root $\alpha \in R$.  Such coweights 
are the concern of this section.

The purpose of
this section is to present an analogue of Proposition~4.4 from \cite{TFSV}
using $\Theta^-_\lambda$ instead of $\Theta_\lambda$. For simplicity, the
theorem is given here for affine Hecke algebras with trivial parameter
systems. The generalization to arbitrary parameter systems is
straightforward (see \cite{TFSV} for notation and details).
Similar arguments to those appearing here apply to $\Theta_\lambda$, giving
a
short proof of Proposition~4.4 from \cite{TFSV}.

\begin{theorem}\label{SameTranslation}
Let $\mu^-$ be minuscule and antidominant, and
$\lambda \in W_0(\mu^-)$. Then
$$
\Theta^-_{\lambda} =
\sum_{x \ : \ \lambda(x) = \lambda}
\widetilde{R}_{x,t_\lambda}(Q) \tilde{T}_x.
$$
\end{theorem}

We begin with some lemmas. For a proof of the first lemma, refer to 
Proposition~3.4 of \cite{TFSV}, where a similar result is given (see also the proof of Corollary \ref{sheaf_prop3.4}).

\begin{lemma}\label{prop3.4}
Let $\mu^-$ be an antidominant and minuscule coweight, and let
$\tau \in \Omega$ be the unique element such that
$t_{\mu^-} \in W_a \tau$. Let $\lambda \in W_0(\mu^-)$. Suppose that
$\lambda - \mu^-$ is a sum of $p$ simple coroots
$(0 \leq p \leq l(t_{\mu^-}) = r)$. Then there exists a sequence of simple
roots
$\alpha_1,\ldots,\alpha_p$ such that the following hold (setting
$s_i = s_{\alpha_i}$):

\noindent {\rm (1)}
$\langle \alpha_i,s_{i-1} \cdots s_1(\mu^-) \rangle = -1,
\forall 1 \leq i \leq p$;

\noindent {\rm (2)}
there is a reduced expression for $t_{\mu^-}$ of the form
$t_{\mu^-} = s_1 \cdots s_p t_1 \cdots t_{r-p} \tau$;

\noindent {\rm (3)}
there is a reduced expression for $t_\lambda$ of the form
$t_\lambda = t_1 \cdots t_{r-p}
(^\tau \! s_1) \cdots (^\tau \! s_p) \tau$;

\noindent {\rm (4)}
$\Theta^-_\lambda =
\tilde{T}_{t_1} \cdots \tilde{T}_{t_{r-p}}
\tilde{T}^{-1}_{^\tau \! s_1} \cdots \tilde{T}^{-1}_{^\tau \! s_p}
\tilde{T}_\tau$,

\noindent where $t_j \in S_a, \forall j \in \{ 1,\ldots,r-p \}$.
\end{lemma}

\begin{lemma}\label{AlwaysAdd}
Let $x \in \widetilde{W}$, and suppose that $xs_{\alpha} > x$ for all
$\alpha \in \Pi$. Then $l(xw) = l(x) + l(w)$ for all $w \in W_0$.
\end{lemma}
Since $\mu^-$ is antidominant, this lemma applies to the expression
$t_{\mu^-} = s_1 \cdots s_p t_1 \cdots t_{r-p} \tau$. It then
applies to the expression $t_1 \cdots t_{r-p} \tau$ as well.
It follows that we can think of the formula in Lemma~\ref{prop3.4}(4) as
$$
\Theta^-_\lambda =
\tilde{T}_{w^\lambda}
\tilde{T}^{-1}_{s_1} \cdots \tilde{T}^{-1}_{s_p} =
\tilde{T}_{w^\lambda} \tilde{T}^{-1}_{w^{-1}}
$$
where $t_\lambda = w^\lambda w$, with $w \in W_0$ and
$w^\lambda$ the minimal length representative for the coset
$t_\lambda W_0$. This observation is helpful towards proving the main
result of this section.

\begin{lemma}\label{bijection}
For $\lambda$, $s_1, \ldots, s_p$, and $t_1, \ldots, t_{r-p}$ as in
Lemma~\ref{prop3.4}, the mapping 
$$
\{ y \ : \ y \leq \ s_1 \cdots s_p \} \ \longrightarrow \
\{ x \ : \ \lambda(x) = \lambda \ {\rm and} \ x \leq t_\lambda \}
$$
defined by
$$
y \ \mapsto \ t_1 \cdots t_{r-p} \tau y = w^\lambda y
$$
is bijective.
\end{lemma}

\medskip

\noindent \begin{proof2}
Let $w = s_1 \cdots s_p$, and $w^\lambda = t_1 \cdots t_{r-p} \tau$, so that
$t_{\mu^-} = w w^\lambda$ and $t_\lambda = w^\lambda w$. We have
$$
\tilde{T}^{-1}_{s_1} \cdots \tilde{T}^{-1}_{s_p} =
\sum_{y \ : \ y \leq w}
\widetilde{R}_{y,w} (Q) \tilde{T}_y
$$
The expression for $\Theta^-_\lambda$ of Lemma~\ref{prop3.4},
together with the fact that
$\tilde{T}_{w^\lambda} \tilde{T}_y = \tilde{T}_{w^\lambda y}$ for all $y \in 
W_0$
(since $l(w^\lambda y) = l(w^\lambda) + l(y)$), implies
$$
\Theta^-_{\lambda} =
\sum_{y \ : \ y \leq w}
\widetilde{R}_{y,w} (Q) \tilde{T}_{w^\lambda y}.
$$
Using the recursion formula of Lemma~2.5 (1) from \cite{TFSV}, we obtain
$\widetilde{R}_{y,w}(Q) = \widetilde{R}_{w^\lambda y, t_\lambda}(Q)$.
In view of
the bijection given in Lemma~\ref{bijection}, we have
$$
\sum_{y \ : \ y \leq w}
\widetilde{R}_{w^\lambda y,t_\lambda} (Q) \tilde{T}_{w^\lambda y} =
\sum_{x \ : \ \lambda(x) = \lambda}
\widetilde{R}_{x,t_\lambda} (Q) \tilde{T}_x,
$$
which completes the proof.
\end{proof2}

\medskip

For the minuscule case, Theorem \ref{SameTranslation} yields the
following improvement on Lemma \ref{support}.  

\begin{corollary}\label{support2}
Let $\lambda \in X_*$ be minuscule.
Then
$$
{\rm supp}(\Theta^-_\lambda) = \{ x \ : \ \lambda(x) = \lambda
\ {\rm and} \ x \leq t_\lambda \}.
$$
\end{corollary}
Here we have used Lemma 2.5 (5) of \cite{TFSV}, which asserts that $\tilde{R}_{x,y}(Q) \neq 0$ if and only if $x \leq y$.

The Bernstein function $z_\mu$ has a very simple form when $\mu$ is 
minuscule (cf. Theorem 4.3 of \cite{TFSV}):

\begin{corollary}\label{bernsteinformula1}
If $\mu$ is dominant and minuscule, then
$$
z_\mu = \sum_{x \in {\rm Adm}(\mu)} \tilde{R}_{x,t_{\lambda(x)}}(Q) 
\tilde{T}_x.
$$
\end{corollary}

\section{Multiples of the Drinfeld case}
Fix positive integers $n$ and $m$, and an integer $1 \leq k \leq n$.
In this section, we establish in Theorem~\ref{Drinfeld} a formula for
the $\Theta^-_\lambda$ functions of ${\rm GL}_n$ when $\lambda = me_k$
(where $e_k$ denotes the coweight of ${\rm GL}_n$ with $k$th
coordinate equal to $1$, and all other coordinates equal to $0$).

In this section, we adopt the following notation: for $1 \leq i \leq n-1$,
let $\alpha_i = \check{\alpha}_i = e_i - e_{i+1}$, and let $s_i = 
s_{\alpha_i}$.  We single out the element $\tau \in \Omega$ given by $\tau =
t_{(1,0,\dots,0)} s_1\cdots s_{n-1}$.

\begin{theorem}\label{Drinfeld}
For the coweight $me_k$ of ${\rm GL}_n$, we have
$$
\Theta^-_{me_k} =
\sum_{x \ : \ \lambda(x) \preceq me_k}
\widetilde{R}_{x,t_{me_k}}(Q) \tilde{T}_x.
$$
\end{theorem}

Consequently, we have a result for $me_k$
analogous to that of Corollary~\ref{support2} for $\lambda$ minuscule,
that is,
$$
{\rm supp}(\Theta^-_{me_k}) = \{ x \ : \ \lambda(x) \preceq me_k \
{\rm and} \ x \leq t_{me_k} \}.
$$
We also get the following explicit formula for the Bernstein function 
$z_{me_1}$, analogous to Corollary \ref{bernsteinformula1}:

\begin{corollary}\label{bernsteinformula2}
$$
z_{me_1} = \sum_{x \in {\rm Adm}(\mu)} \left( \sum_{\lambda(x) \preceq 
\lambda \,\, , \,\, x \leq t_\lambda}
\tilde{R}_{x,t_\lambda}(Q) \right) \tilde{T}_x,
$$
where the inner sum ranges over $\lambda \in W_0(me_1)$ such that
$\lambda(x) \preceq \lambda$ and $x \leq t_\lambda$.
\end{corollary}

\medskip

We require three lemmas before the proof of the theorem.  In the following arguments we use the notation $\prod$ to denote products even though we are working in a non-commutative ring.  We will use the following convention: $\prod^n_{i=1} a_i$ will denote the product 
$a_1a_2 \cdots a_n$ (in that order).

\begin{lemma}\label{theta}
For the coweight $me_k$ of ${\rm GL}_n$, we have
$$
\Theta^-_{me_k} =
(\tilde{T}_{s_{k-1}} \cdots \tilde{T}_{s_1} \tilde{T}_\tau
\tilde{T}^{-1}_{s_{n-1}} \cdots \tilde{T}^{-1}_{s_k})^m.
$$
\end{lemma}
\begin{proof}
From Lemma~\ref{prop3.4}, we have
$$
\Theta^-_{e_k}
= \tilde{T}_{s_{k-1}} \cdots \tilde{T}_{s_1}
\tilde{T}_\tau \tilde{T}^{-1}_{s_{n-1}} \cdots \tilde{T}^{-1}_{s_k}.
$$
Then the formula for
$\Theta^-_{me_k}$ follows from the fact that
$\Theta^-_{\nu_1 + \nu_2} = \Theta^-_{\nu_1} \Theta^-_{\nu_2}$ for all
$\nu_1,\nu_2~\in~X_*$.
\end{proof}

\begin{lemma}\label{collapse}
Let $w, y \in \tilde{W}$. Then
$$
\tilde{T}_w \tilde{T}_y =
\sum_{x = \widetilde{w} \widetilde{y}} a_x \tilde{T}_x
$$
where $\widetilde{w}$ and $\widetilde{y}$
range over certain subexpressions of $w$ and $y$,
respectively.
\end{lemma}
\begin{proof}
This is an easy induction on the length of $y$.
\end{proof}

\begin{lemma}\label{specialExpressions}
Let x be a subexpression of
$$
t_{me_k} = (s_{k-1} \cdots s_1 \tau s_{n-1} \cdots s_k)^m
$$
such that for some $1 \leq i \leq k-1$,
at least one $s_i$ is deleted.
Then $\lambda(x) \npreceq me_k$.
\end{lemma}
\begin{proof}
We can write
$$
x = u_1 \tau v_1 \cdots u_m \tau v_m
$$
for suitable subexpressions $u_1, \ldots, u_m$ and $v_1,\ldots,v_m$
of $s_{k-1} \cdots s_1$ and $s_{n-1} \cdots s_k$, respectively.
Suppose that $p$ is the
least index such that $u_p \not= s_{k-1} \cdots s_1$.
Then
$$
x = \bigg( \prod_{i=1}^{p-1} t_{e_k} s_k \cdots s_{n-1} v_i \bigg)
(t_{e_j} u_p s_1 \cdots s_{n-1} v_p)
\bigg( \prod_{i=p+1}^m u_i \tau v_i \bigg),
$$
for some $j$ with $1 \leq j < k$.
Since $s_k \cdots s_{n-1} v(e_j) = e_j$ for any subexpression
$v$ of $s_{n-1} \cdots s_k$,
$$
\bigg( \prod_{i=1}^{p-1} t_{e_k} s_k \cdots s_{n-1} v_i \bigg) (t_{e_j}) \bigg( \prod_{i=1}^{p-1} t_{e_k} s_k \cdots s_{n-1} v_i \bigg)^{-1}
= t_{e_j}.
$$
It follows that the translation part $\lambda(x)$ is the sum of $e_j$ and a non-negative integral linear combination of vectors $e_i$ ($i=1,\dots,n$).  Indeed, the translation part of
$$\bigg( \prod_{i=1}^{p-1} t_{e_k} s_k \cdots s_{n-1} v_i \bigg) (u_p s_1 \cdots s_{n-1} v_p)
\bigg( \prod_{i=p+1}^m u_i \tau v_i \bigg)
$$
is necessarily a vector $(b_1,\dots,b_n)$ where $b_i \in {\mathbb Z}_+$ for every $i$.

We thus see that one of the first $k-1$ coordinates of
$\lambda(x)$ is positive (namely the $j$-th coordinate is), and this implies that $\lambda(x) \npreceq me_k$.
\end{proof}

\noindent \begin{proof3}
Let
$E = (\{ 0,1 \}^{n-1})^m - ( \{ 1 \}^{k-1} \times \{ 0,1 \}^{n-k} )^m$.
Using $\tilde{T}^{-1}_s = \tilde{T}_s + Q$ and expanding the left hand side,
we can write
\begin{eqnarray}\label{thetaM}
\tilde{T}^{-1}_{-me_k} =
\Theta^-_{me_k} +
\sum_{\epsilon \in E} Q^{m(n-1) -\sigma(\epsilon)}
\prod_{i=1}^{m}
\tilde{T}^{\epsilon^i_{k-1}}_{s_{k-1}} \cdots \tilde{T}^{\epsilon^i_1}_{s_1}
\tilde{T}_\tau
\tilde{T}^{\epsilon^i_{n-1}}_{s_{n-1}} \cdots 
\tilde{T}^{\epsilon^i_k}_{s_k}.
\end{eqnarray}
Here $\epsilon^i_j \in \{ 0,1 \}$ and we denote
$\epsilon = (\epsilon^1, \ldots, \epsilon^m)$, and
$\epsilon^i = (\epsilon^i_1, \ldots, \epsilon^i_{n-1})$,
$i = 1,\ldots,m$,
and $\sigma(\epsilon) = \sum_{i=1}^m \sum_{j=1}^{n-1} \epsilon^i_j$.

From Lemma~\ref{support}, we know that
supp$(\Theta^-_{me_k}) \subset \{ x \ : \ \lambda(x) \preceq me_k \}$. Thus
we need only prove that if $\tilde{T}_x$ is in the support of the second
term on the right hand side, then
$\lambda(x) \npreceq me_k$.  Indeed, then the first and second terms on the right hand side of (\ref{thetaM}) have disjoint supports, and so the coefficients of like terms will be equal in $\tilde{T}^{-1}_{-me_k}$ and $\Theta^-_{me_k}$.

Let $\epsilon = (\epsilon^1, \ldots, \epsilon^m) \in E$, and consider
$$
\prod_{i=1}^{m}
\tilde{T}^{\epsilon^i_{k-1}}_{s_{k-1}} \cdots \tilde{T}^{\epsilon^i_1}_{s_1}
\tilde{T}_\tau
\tilde{T}^{\epsilon^i_{n-1}}_{s_{n-1}} \cdots
\tilde{T}^{\epsilon^i_k}_{s_k}.
$$
By Lemma~\ref{collapse}, if $x$ is in the support of this
product, then $x$ is a subexpression of
$$
\prod_{i=1}^m
s_{k-1}^{\epsilon^i_{k-1}} \cdots s_1^{\epsilon^i_1}
\tau
s_{n-1}^{\epsilon^i_{n-1}} \cdots s_k^{\epsilon^i_k}.
$$
Since $E$ excludes the elements of
$(\{ 1 \}^{k-1} \times \{ 0,1 \}^{n-k} )^m$, we know that for
some $1 \leq i \leq m$ and $1 \leq j \leq k-1$, that
$\epsilon^i_j = 0$. But this is equivalent to the deletion of some $s_j$ ($1 
\leq j \leq k-1$) from
the
expression $t_{me_k}~=~(s_{k-1} \cdots s_1 \tau s_{n-1} \cdots s_k)^m$.
By Lemma~\ref{specialExpressions}, we have $\lambda(x) \npreceq me_k$, and 
the proof is complete.
\end{proof3}

\section{Minimal expressions}

We say $\Theta^-_\lambda$ has a {\em minimal expression} if it can be written in the form
$$
\Theta^-_\lambda = \tilde{T}^{\epsilon_1}_{t_1} \cdots \tilde{T}^{\epsilon_r}_{t_r}\tilde{T}_\tau,
$$
where $t_\lambda = t_1 \cdots t_r \tau$ ($t_i \in S_a, \, \tau \in \Omega$) is a reduced expression and $\epsilon_i \in \{ 1,-1 \}$ for every $1 \leq i \leq r$.  Such expressions played a key role in Theorems 3.1 and 4.1.

Lemma 3.2 asserts that $\Theta^-_\lambda$ has a minimal expression whenever $\lambda$ is minuscule.  If $\lambda$ is any coweight for ${\rm GL}_n$, then we may write
$$
\lambda = \lambda_1 + \cdots + \lambda_k
$$
where each $\lambda_j$ is minuscule and 
$$
l(t_\lambda) = l(t_{\lambda_1}) + \cdots + l(t_{\lambda_k}).
$$
It follows that for any coweight $\lambda$ of ${\rm GL}_n$, there is a minimal expression for $\Theta^-_\lambda$.  Letting $w^{\lambda_i}$ denote the minimal representative for the coset $t_{\lambda_i}W_0$ and writing $t_{\lambda_i} = w^{\lambda_i}w_i$ ($w_i \in W_0$), we may recover a minimal expression for
$$
\Theta^-_\lambda = \tilde{T}_{w^{\lambda_1}}\tilde{T}^{-1}_{w^{-1}_1} \cdots
\tilde{T}_{w^{\lambda_k}}\tilde{T}^{-1}_{w^{-1}_k}
$$
by choosing reduced expressions for every $w^{\lambda_i}$ and $w_i$.

Clearly a similar result would follow for any root system with the property that $\Theta^-_{\omega}$ has a minimal expression for every Weyl conjugate $\omega$ 
of every fundamental coweight.  It seems to be an interesting combinatorial problem 
to determine the root systems (besides that for ${\rm GL}_n$) which satisfy this property. 

In principle, a minimal expression for $\Theta^-_\lambda$ allows one to write it as an explicit linear combination of the Iwahori-Matsumoto generators $T_w$, simply by using the formula $\tilde{T}^{-1}_s = \tilde{T}_s + Q$ and expanding the product.  The result is a linear combination of certain products
$$
T_{s_1} \cdots T_{s_g} T_\sigma
$$
($s_i \in S_a, \, \sigma \in \Omega$), where $s_1 \cdots s_g \sigma$ 
ranges over certain subexpressions of $t_\lambda$ (which subexpressions occur is governed by the signs $\epsilon_j$ in the minimal expression). These may in turn be simplified by using the well-known formula
$$
T_{s_1} \cdots T_{s_g} = \sum_w N(\underline{s},w,q) T_w,
$$
(cf. \cite{Lus1}, Lemma 3.7).  Here $\underline{s} = (s_1,\dots, s_g)$ and 
$N(\underline{s},w,q)$ is the number of ${\mathbb F}_q$-rational points on the variety $Z(\underline{s},w)$ consisting of all sequences $(I_1,\dots,I_g)$ where the $I_i$ are Iwahori subgroups of $G^{sc}(\bar{{\mathbb F}}_q(\!(t)\!))$ (here $G^{sc}$ is the simply connected group associated to the given root system) such that the relative positions of adjacent subgroups satisfy
$$
{\rm inv}(I_{i-1},I_i) = s_i
$$
for all $1 \leq i \leq g$, and $I_g = wI_0w^{-1}$, where $I_0$ is a fixed ``standard'' Iwahori subgroup.  

We forgo the cumbersome task of describing more completely the resulting expressions for $\Theta^-_\lambda$ in terms of the generators $T_w$.  The combinatorics are best described in the geometric framework of Demazure resolutions.  We explain this in the following section.

\medskip

\noindent {\em Remark.}  
Let $\lambda$ be a coweight for ${\rm GL}_n$, and write $\lambda = m_1e_1 + \cdots + m_ne_n$.  One finds a similar 
expression for $\Theta^-_\lambda$ starting from
$$
\Theta^-_\lambda = \Theta^-_{m_1e_1} \cdots \Theta^-_{m_ne_n},
$$
and making use of Theorem 4.1.

\section{Sheaf-theoretic meaning of minimal expressions}

The goal of this section is to describe a sheaf-theoretic interpretation of a minimal expression for $\Theta^-_\lambda$: the corresponding perverse sheaf on the affine flag variety is the push-forward of an explicit perverse sheaf on a Demazure resolution of the Schubert variety $X(t_\lambda)$.  We proceed to illustrate this statement in more detail.

\subsection{Affine flag variety}

Let $k = {\mathbb F}_q$ denote the finite field with $q$ elements, and let $\bar{k}$ denote an algebraic closure of $k$.  Let $G$ be the split connected reductive group over $k$ whose root system is $(X^*,X_*,R, \check{R}, \Pi)$.  Choose a split torus $T$ and a $k$-rational Borel subgroup $B$ containing $T$, which give rise to $R$ and $\Pi$.  

Denote by ${\mathcal Fl}$ the affine flag variety for $G$.  This is an ind-scheme over $k$ whose $k$-points are given by
$$
{\mathcal Fl}(k) = G(k(\!(t)\!))/I_{k},
$$
where $I = I_k \subset G(k[[t]])$ is the Iwahori subgroup whose reduction 
modulo $t$ is $B$.

Fix a prime $\ell \neq \mbox{char} (k)$, and make a fixed choice for $\sqrt{q} \in \bar{\mathbb Q}_\ell$ (for Tate twists).

Let $D^b({\mathcal Fl})$ denote the category $D^b_c({\mathcal Fl}, \bar{\mathbb Q}_\ell)$.  By definition $D^b_c({\mathcal Fl},\bar{\mathbb Q}_\ell)$ is the inductive 2-limit of categories $D^b_c(X, \bar{\mathbb Q}_\ell)$ where $X \subset {\mathcal Fl}$ ranges over all projective $k$-schemes which are closed subfunctors of the ind-scheme ${\mathcal Fl}$.  The category $D^b_c(X,{\mathbb Q}_\ell)$ is the ``derived'' category of Deligne \cite{Weil2}: ${\mathbb Q}_\ell \otimes$ the projective 2-limit of the categories $D^b_{ctf}(X, {\mathbb Z} / \ell^n {\mathbb Z})$.  For any finite extension $E$ of ${\mathbb Q}_\ell$ contained in $\bar{\mathbb Q}_\ell$, the definition of $D^b_c(X,E)$ is similar, and by definition $D^b_c(X,\bar{\mathbb Q}_\ell)$ is the inductive 2-limit of the categories $D^b_c(X,E)$.

For $f: X \rightarrow Y$ a morphism of finite-type $k$-schemes, we have the four ``derived'' functors $f_*, f_! : 
D^b_c(X,\bar{\mathbb Q}_\ell) \rightarrow D^b_c(Y,\bar{\mathbb Q}_\ell)$ and $f^*, f^! : D^b_c(Y,\bar{\mathbb Q}_\ell) \rightarrow D^b_c(X,\bar{\mathbb Q}_\ell)$.  This notation should cause no confusion, since we never use the non-derived versions of the pull-back and push-forward functors in this paper.

We define the category $P_I({\mathcal Fl})$: it is the full subcategory of 
$D^b({\mathcal Fl})$ whose objects are $I$-equivariant perverse sheaves for the middle perversity (by definition the latter have finite dimensional support). 

The $I$-orbits on ${\mathcal Fl}$ correspond to $\widetilde{W}$.  Given $w \in \widetilde{W}$, we denote by $Y(w) = IwI/I$ the corresponding Bruhat cell, and we denote its closure by $X(w) = \overline{Y(w)}$.  Further, let $\bar{{\mathbb Q}}_{\ell,w}$ denote the constant sheaf on $Y(w)$, and 
define ${\mathcal A}_w = \bar{{\mathbb Q}}_{\ell,w}[l(w)](l(w)/2)$.
This is a self-dual perverse sheaf on $Y(w)$.  

Let $j_w: Y(w) \hookrightarrow X(w)$ denote the open immersion.  We define
$J_{w*} = j_{w*}{\mathcal A}_w$ and $J_{w!} = j_{w!}{\mathcal A}_w$.   These are perverse sheaves in $P_I({\mathcal Fl})$ satisfying $D(J_{w*}) = J_{w!}$.  
(Here $D$ denotes Verdier duality.)

Given ${\mathcal G} \in P_I({\mathcal Fl})$ we may define the corresponding function 
$[{\mathcal G}]$ on ${\mathcal Fl}(k)$, which we may identify with an element in 
${\mathcal H}$:
$$
[{\mathcal G}](x) = {\rm Tr}({\rm Fr}_q, {\mathcal G}_x),
$$
where ${\rm Fr}_q$ denotes the Frobenius morphism on ${\mathcal Fl}_{\bar{k}}$ (raising coordinates to power $q$).  We have 
\begin{align*}
[J_{w!}]  &= \varepsilon_w q_w^{-1/2} T_w \\
[J_{w*}]  &= \varepsilon_w q_w^{1/2} T^{-1}_{w^{-1}}.
\end{align*}

\subsection{Convolution of sheaves}

Following Lusztig \cite{Lus3}, one can define a convolution product 
$\star \, : P_I({\mathcal Fl}) \times P_I({\mathcal Fl}) \rightarrow D^b({\mathcal Fl})$.  We formulate this in a way similar to \cite{Haines-Ngo}.  Given ${\mathcal G}_i \in P_I({\mathcal Fl})$, $i=1,2$, we can choose $X(w_i)$ such that 
the support of ${\mathcal G}_i$ is contained in $X(w_i)$, 
for $i=1,2$.  We may identify ${\mathcal Fl}$ with the space of all ``affine flags'' ${\mathcal L}$ for $G(k(\!(t)\!))$; there is a base point ${\mathcal L}_0$ whose stabilizer in $G(k(\!(t)\!))$ is the ``standard'' Iwahori subgroup $I$.  
Then $X(w)$ is identified with the space of all affine flags ${\mathcal L}$ such that the relative position between the base point ${\mathcal L}_0$ and 
${\mathcal L}$ satisfies 
$$
{\rm inv}({\mathcal L}_0,{\mathcal L}) \leq w
$$
in the Bruhat order on $\widetilde{W}$.  The ``twisted'' product $X(w_1) 
\tilde{\times} X(w_2)$ is the space of pairs $({\mathcal L},{\mathcal L}') \in {\mathcal Fl} \times {\mathcal Fl}$ such that
\begin{align*}
{\rm inv}({\mathcal L}_0,{\mathcal L}) &\leq w_1  \\
{\rm inv}({\mathcal L},{\mathcal L}') &\leq w_2.
\end{align*}
We can find a finite-dimensional projective subvariety $X \subset {\mathcal Fl}$ with the property that 
$({\mathcal L},{\mathcal L}') \in X(w_1) \tilde{\times} X(w_2) \Rightarrow {\mathcal L}' \in X$.
The ``multiplication'' map $m: X(w_1) \tilde{\times} X(w_2) \rightarrow X$ given by $({\mathcal L},{\mathcal L}') \mapsto {\mathcal L}'$ is proper.  

Now ${\mathcal G}_i$ ($i=1,2$) determine a well-defined perverse sheaf 
${\mathcal G}_1 \tilde{\boxtimes} {\mathcal G}_2$ on $X(w_1) \tilde{\times} X(w_2)$ (see e.g. \cite{Haines-Ngo}).  We define
$$
{\mathcal G}_1 \star {\mathcal G}_2 = m_*({\mathcal G}_1 \tilde{\boxtimes} 
{\mathcal G}_2).
$$
The convolution ${\mathcal G}_1 \star {\mathcal G}_2 \in D^b({\mathcal Fl})$ is independent of the choice of $X(w_i)$ and $X$.

The object ${\mathcal G}_1 \star {\mathcal G}_2$ is $I$-equivariant in a suitable sense, so that we can regard its function ${\rm Tr}({\rm Fr}_q , {\mathcal G}_1 \star{\mathcal G}_2)$ as an element of the Hecke algebra ${\mathcal H}$.

It is well-known that this product is compatible with the function-sheaf dictionary:
$$
[{\mathcal G}_1 \star {\mathcal G}_2] = [{\mathcal G}_1] \star [{\mathcal G}_2].
$$
Here $\star$ on the right hand side is just the usual product in ${\mathcal H}$.

Later we shall use the following fact, referred to in the sequel simply as {\em associativity}: if ${\mathcal G}_i$ ($i = 1,2,3$) are objects of $P_I({\mathcal Fl})$ such that ${\mathcal G}_1 \star {\mathcal G}_2 \in 
P_I({\mathcal Fl})$ and ${\mathcal G}_2 \star {\mathcal G}_3 \in P_I({\mathcal Fl})$, then there is a canonical isomorphism ${\mathcal G}_1 \star ({\mathcal G}_2 \star {\mathcal G}_3) 
\tilde{\rightarrow} ({\mathcal G}_1 \star {\mathcal G}_2) \star {\mathcal G}_3$ (the ``associativity constraint'').  This is proved by identifying each canonically with the ``triple product'' 
${\mathcal G}_1 \star {\mathcal G}_2 \star {\mathcal G}_3$, whose construction is similar (see section 6.3).

\subsection{Demazure resolution vs. twisted product}

It is clear that we can define in a similar way the $k$-fold convolution product
$\star \, : P_I({\mathcal Fl})^k \rightarrow D^b({\mathcal Fl})$.  To do this we define the $k$-fold
twisted product $X(w_1) \tilde{\times} \cdots \tilde{\times} X(w_k)$ to be the space of $k$-tuples $({\mathcal L}_1,\dots,{\mathcal L}_k)$ such that 
$$
{\rm inv}({\mathcal L}_{i-1},{\mathcal L}_i) \leq w_i
$$
for $1 \leq i \leq k$.  If $w = w_1\cdots w_k$ and 
$l(w) = \sum_i l(w_i)$, then the multiplication map
$$
m: X(w_1) \tilde{\times} \cdots \tilde{\times} X(w_k) \longrightarrow X(w)
$$
given by $({\mathcal L}_1,\dots,{\mathcal L}_k) \mapsto {\mathcal L}_k$ induces an isomorphism on open subschemes
$$
m: Y(w_1) \tilde{\times} \cdots \tilde{\times} Y(w_k) \tilde{\longrightarrow} Y(w).
$$
If moreover each $w_i$ is a simple reflection $s_i$, then the twisted product is smooth, since it is a succession of ${\mathbb P}^1$-bundles.  We have proved the following lemma.

\begin{lemma}\label{demazure}  If $w = s_1 \cdots s_r$ is a reduced expression, then the twisted product$$
m: X(s_1) \tilde{\times} \cdots \tilde{\times} X(s_r) \longrightarrow X(w)
$$
is a Demazure resolution for the Schubert variety $X(w)$.  
\end{lemma}

\medskip

We use the $k$-fold twisted product to define the $k$-fold convolution product: as before, the objects ${\mathcal G}_i \in P_I({\mathcal Fl})$ determine a (unique, perverse) twisted exterior product ${\mathcal G}_1 \tilde{\boxtimes} \cdots \tilde{\boxtimes} {\mathcal G}_k$, and we set
$$
{\mathcal G}_1 \star \cdots \star {\mathcal G}_k = 
m_*({\mathcal G}_1 \tilde{\boxtimes} \cdots \tilde{\boxtimes} {\mathcal G}_k).
$$

We have the following generalized associativity constraint.  Consider a product with $k$ terms
$$
(\cdots( {\mathcal G}_{1} \star (\cdots \star {\mathcal G}_{i}))\cdots) \star (\cdots( {\mathcal G}_{j} \star (\cdots \star {\mathcal G}_{k}))\cdots)
$$
where the placement of the parentheses is arbitrary with the proviso that the product is defined (i.e., at every stage we convolve objects of $P_I({\mathcal Fl})$).  Then this can be identified canonically with the $k$-fold product
$$
{\mathcal G}_1 \star \cdots \star {\mathcal G}_k.
$$
This can be seen easily by induction on $k$.

\subsection{Properties of certain convolutions}

The convolution of two $I$-equivariant perverse sheaves on ${\mathcal Fl}$ is not perverse in general.  However, the following result of I. Mirkovic 
(unpublished) shows this conclusion does hold in some important cases.
We are grateful to R. Bezrukavnikov, who communicated this result to the first author, and to I. Mirkovic, for his kind permission to include the result in this paper.  

In the notation of \cite{BBD}, we let $ ^pD^{\leq 0}({\mathcal Fl})$ (resp. 
$ ^pD^{\geq 0}({\mathcal Fl})$) denote the objects $P \in D^b({\mathcal Fl})$ whose perverse cohomology sheaves vanish in degree $\geq 1$ (resp. $\leq -1$).  Thus the perverse sheaves on ${\mathcal Fl}$ are precisely the objects in 
$^pD^{\geq 0}({\mathcal Fl}) \, \cap \, ^pD^{\leq 0}({\mathcal Fl})$. 

\begin{proposition}[Mirkovic] \label{Mirkovic}
\noindent {\em (a)}  Let $P \in P_I({\mathcal Fl})$. Then
for any $w \in \widetilde{W}$, we have 
\begin{enumerate}
\item[(1)] $J_{w!} \star P$ and $P \star J_{w!}$ belong to $ ^pD^{\geq 0}({\mathcal Fl})$,
\item[(2)] $J_{w*} \star P$ and $P \star J_{w*}$ belong to $ ^pD^{\leq 0}({\mathcal 
Fl})$.
\end{enumerate} 

\noindent {\em (b)}  In particular, $J_{w_1*} \star J_{w_2!}$ and $J_{w_1!} \star J_{w_2*}$ are perverse, for every 
$w_1, w_2 \in \widetilde{W}$.  

\smallskip

\end{proposition}

\begin{proof}  
\noindent (a),(1).  We consider $J_{w!} \star P$.  Suppose $P$ is supported on $X(w')$ and recall that the convolution is given by 
$$
J_{w!} \star P = m_!(J_{w!} \tilde{\boxtimes} P),
$$
where $m: X(w) \tilde{\times} X(w') \rightarrow X$ is as in section 6.2 (note that 
$m_! = m_*$, since $m$ is proper).
We have
\begin{align*}
J_{w!} \star P & = m_!(j_{w!}({\mathcal A}_w) \tilde{\boxtimes} P) \\ 
&= [m \circ (j_{w} \tilde{\times} {\rm id})]_!({\mathcal A}_w \tilde{\boxtimes} P).
\end{align*}
Note that $m \circ (j_{w} \tilde{\times} {\rm id}) : 
Y(w) \tilde{\times} X(w') \rightarrow X$ is affine, and that ${\mathcal A}_w \tilde{\boxtimes} P$ is perverse on its source.  Therefore 
$J_{w!} \star P \in \, ^pD^{\geq 0}({\mathcal Fl})$ follows from the following general fact (cf. \cite{BBD}, Thm. 4.1.1 and Cor. 4.1.2): If $F:X \rightarrow Y$ is an affine morphism, then 
$F_*: D^b_c(X) \rightarrow D^b_c(Y)$ (resp. $F_!$) preserves $^pD^{\leq 0}$ (resp. $^pD^{\geq 0}$).  

A similar argument gives $P \star J_{w!} \in \, ^pD^{\geq 0}({\mathcal Fl})$, once it is noted that 
$m \circ ({\rm id} \tilde{\times} j_{w}) :  X(w') \tilde{\times} Y(w) \rightarrow X$ is also an affine morphism.  

Part (a),(2) is similar, and part (b) is an immediate consequence of part (a).

\end{proof}

Now we let $P_I({\mathcal Fl}) \cap J_{w*}P_I({\mathcal Fl})$ 
(resp. $P_I({\mathcal Fl}) \cap J_{w!}P_I({\mathcal Fl})$) denote the full subcategory of $P_I({\mathcal Fl})$ whose objects are of the form $J_{w*}\star{\mathcal P}$ (resp. $J_{w!}\star{\mathcal P}$) for some 
${\mathcal P} \in P_I({\mathcal Fl})$. 

\begin{corollary}\label{equivalence} For any $w \in \widetilde{W}$, we have 
$$
J_{w!} \star J_{w^{-1}*} = J_{w^{-1}*} \star J_{w!} = J_e,
$$
where $e \in \widetilde{W}$ is the identity element.
Together with the associativity constraint, these identities imply that
$$
J_{w!} \star - \, : \, P_I({\mathcal Fl}) \cap J_{w^{-1}*}P_I({\mathcal Fl}) \longrightarrow P_I({\mathcal Fl}) \cap J_{w!}P_I({\mathcal Fl})
$$
is an equivalence of categories, with inverse $J_{w^{-1}*} \star -$.
\end{corollary}

\begin{proof}
We prove that $J_{w!} \star J_{w^{-1}*} = J_e$ (the other equality is similar).  Let $X(y)$ be an irreducible component in the support of 
$P := J_{w!} \star J_{w^{-1}*}$.  Since $P$ is perverse and $I$-equivariant, the restriction $P|Y(y)$ is an $I$-equivariant $\ell$-adic local system on the affine space $Y(y)$.  We need to show that $y = e$ and 
$P|Y(e) = \bar{\mathbb Q}_\ell$.  For the former it is sufficient to prove that $y \neq e$ implies $P|Y(y) = 0$.

Since $P|Y(y)$ is an $I$-equivariant $\ell$-adic local system on the $I$-orbit $Y(y)$, and the stabilizer in $I$ of any point in this orbit is geometrically connected, it follows that $P|Y(y)$ is a {\em constant} local system (placed in degree $-l(y)$ when regarded as a complex).  Write $\alpha_1,\alpha_2, \dots, \alpha_r$ for the eigenvalues of ${\rm Fr}_q$ on $P|Y(y)[-l(y)]$, counted with multiplicity.  We have the 
following identity for every $n \geq 1$:
$$
J_e(y) = {\rm Tr}({\rm Fr}_{q^n}, J_{w!}) \star 
{\rm Tr}({\rm Fr}_{q^n}, J_{w^{-1}*})(y) = {\rm Tr}({\rm Fr}_{q^n},P|Y(y)) = \epsilon_y \sum_{i=1}^r \alpha^n_i.
$$ 
We thus have, for every $n \geq 1$: 
\begin{equation*}
\epsilon_y \sum_{i=1}^r \alpha^n_i = 
\begin{cases}
0, \,\,\,\, \mbox{if $y \neq e$} \\
1,  \,\,\,\, \mbox{if $y =e$}.
\end{cases}
\end{equation*}
The linear independence of characters (or rather its proof) implies that distinct numbers $\beta \in \bar{\mathbb Q}^{\times}_\ell$ determine linearly independent characters $n \mapsto \beta^n$ on the semi-group of positive integers.  Together with the above formula this is enough information to determine the eigenvalues of ${\rm Fr}_q$ on $P|(Y(y)$ (if $y \neq e$ then $r=0$, and if $y=e$ then $r=1$ and $\alpha_1 = 1$).  Thus 
$P|Y(y) = 0$ if $y \neq e$ and $P|Y(e) = \bar{\mathbb Q}_\ell$, as desired.

\end{proof}

It is straightforward to check the following properties.  
From now on we omit the convolution sign $\star$ in the product of perverse sheaves.

\begin{lemma}\label{properties}
\begin{enumerate}
\item If $l(xy) = l(x) + l(y)$, then $J_{x!}J_{y!} = J_{xy!}$. 
\item Under the same assumption, $J_{x*}J_{y*} = J_{xy*}$.
\end{enumerate}
\end{lemma}

\begin{proof}
The first can be checked from the definitions, and the second follows on applying Verdier duality.
We have used that Verdier duality is compatible with convolution: 
$D({\mathcal G}_1 \star {\mathcal G}_2) = D{\mathcal G}_1 \star D{\mathcal G}_2$, if ${\mathcal G}_i \in P_I({\mathcal Fl})$ $\, (i = 1,2)$. 

\end{proof}

\smallskip

\noindent {\em Remark.} 
We remark that Corollary \ref{equivalence} and Lemma \ref{properties} allow us to perform algebraic manipulations involving the perverse sheaf $J_{w!}$ (resp. $J_{w*}$): essentially it behaves just like its function $\varepsilon_w \tilde{T}_w$ 
(resp. $\varepsilon_w \tilde{T}^{-1}_{w^{-1}}$) (but when multiplying sheaves, one has to take care that they are each perverse).  
For example, we have the following cancellation property.  Let $P_i \in P_I({\mathcal Fl})$ ($i=1,2$) be such that $J_{w*} \star P_i \in P_I({\mathcal Fl})$ ($i=1,2$); then $J_{w*}P_1 \cong J_{w*}P_2$ implies $P_1 \cong P_2$ (multiply both sides by $J_{w^{-1}!}$ and use associativity).  We shall use this several times in the proof of Lemma \ref{moreproperties} below.

\bigskip

\subsection{Sheaf analogue $\Xi^-_\lambda$ of $\Theta^-_\lambda$}

We now define the sheaf-analogue of $\Theta^-_\lambda$.   We write $J_{\lambda*}$ (resp. $J_{\lambda!}$) in place of $J_{t_\lambda*}$ (resp. $J_{t_\lambda!}$).
If $\lambda = \lambda_1 - \lambda_2$, where $\lambda_i$ is anti-dominant ($i=1,2$), then we define
$$
\Xi^-_\lambda = J_{\lambda_1!}J_{-\lambda_2*}.
$$
By Proposition \ref{Mirkovic},  $\Xi^-_\lambda$ is an object of 
$P_I({\mathcal Fl})$.  Moreover it clearly satisfies 
$$
[\Xi^-_\lambda] = \varepsilon_\lambda \Theta^-_\lambda.
$$
By \cite{Bernstein}, we know that ${\rm supp}(\Theta^-_\lambda) \subset \{ x \, | \, x \leq t_\lambda \}$, and hence
that $\Xi^-_\lambda$ is supported on $X(t_\lambda)$.
\medskip

Our next goal is to prove the sheaf-theoretic analogues of the relations in the Bernstein presentation of ${\mathcal H}$, in the following lemma.  
The proof follows Lemma 4.4 of \cite{Lus4} very closely, 
taking into account Proposition \ref{Mirkovic}, Corollary \ref{equivalence} and Lemma \ref{properties}.  
Since a little extra care must be taken in the present context of perverse sheaves, we give detailed arguments for the convenience of the reader.

\begin{lemma}\label{moreproperties}
\begin{enumerate}
\item[(1)] If $x,y \in \widetilde{W}$ commute and $l(xy) = l(x) + l(y)$, then 
$J_{x!}J_{y^{-1}*} = J_{y^{-1}*}J_{x!}$.  In particular,
if $\mu, \lambda \in X_*$ are both dominant or antidominant, then 
$J_{\mu!}J_{-\lambda*} = J_{-\lambda*}J_{\mu!}$.
\item[(2)] $\Xi^-_\lambda$ is independent of the choice of $\lambda_i$ {\em (}
$i=1,2$ {\em )}.
\item[(3)] $\Xi^-_\lambda \Xi^-_\mu = \Xi^-_{\lambda + \mu}$ for $\lambda, \mu \in X_*$.
\item[(4)] If $s=s_\alpha$ {\em (}$\alpha \in \Pi${\em )} and $\langle \alpha, \lambda \rangle = 0$, then 
$$
J_{s*}\Xi^-_\lambda = \Xi^-_\lambda J_{s*}.
$$
Moreover, this object belongs to $P_I({\mathcal Fl})$.
\item[(5)] If $\langle \alpha, \lambda \rangle = -1$, then
$$
J_{s*} \Xi^-_\lambda J_{s*} = \Xi^-_{s\lambda}.
$$
\end{enumerate}
\end{lemma}

\begin{proof}
\noindent (1).  By Lemma \ref{properties}, we have $J_{x!}J_{y!} = 
J_{xy!} = J_{yx!} = J_{y!}J_{x!}$.   The result follows by two applications of Corollary \ref{equivalence}: multiply first on the left and then on the right by $J_{y^{-1}*}$.

\smallskip

\noindent (2).  Let $\lambda = \lambda_1 - \lambda_2 = \lambda'_1 - \lambda'_2$, where $\lambda_i, \lambda'_i$ are antidominant ($i=1,2$).  Then $\lambda'_2 + \lambda_1 = \lambda'_1 + \lambda_2$, and so $J_{\lambda'_2!}J_{\lambda_1!} = J_{\lambda'_1!}J_{\lambda_2!}$.  Arguing as in (1), this yields $J_{\lambda_1!}J_{-\lambda_2*} = 
J_{-\lambda'_2*}J_{\lambda'_1!}$.  Then (1) yields 
$J_{\lambda_1!}J_{-\lambda_2*} = J_{\lambda'_1!}J_{-\lambda'_2*}$.

\smallskip

\noindent (3).  Write $\lambda = \lambda_1 - \lambda_2$ and $\mu = \mu_1 - \mu_2$, where $\lambda_i, \mu_i$ are antidominant ($i=1,2$).  Then by associativity we have 
\begin{align*}
\Xi^-_{\lambda} \Xi^-_{\mu} &= (J_{\lambda_1!}J_{-\lambda_2*})(J_{\mu_1!}J_{-\mu_2*}) \\
&= J_{\lambda_1!} (J_{-\lambda_2*}J_{\mu_1!})J_{-\mu_2*} \\
&= J_{\lambda_1!} (J_{\mu_1!}J_{-\lambda_2*})J_{-\mu_2*} \\
&= (J_{\lambda_1!}J_{\mu_1!})(J_{-\lambda_2*}J_{-\mu_2*}) \\
&= J_{(\lambda_1 + \mu_1)!}J_{-(\lambda_2 + \mu_2)*} \\
&= \Xi^-_{\lambda + \mu}.
\end{align*}

\smallskip

\noindent (4).  We may write $\lambda = \lambda_1 - \lambda_2$, where $\lambda_i$ is antidominant and $\langle \alpha, \lambda_i \rangle = 0$, for $i=1,2$.  Since $s$ commutes with $t_{\lambda_i}$ and $l(st_{\lambda_i}) = l(t_{\lambda_i}s) = l(t_{\lambda_i}) + 1$, the result follows from (1), Lemma \ref{properties}, and associativity.   

We note that $J_{s*}(J_{\lambda_1!}J_{-\lambda_2*}) = J_{st_{-\lambda_2}*}J_{t_{\lambda_1}!}$ is $I$-equivariant and perverse, by Proposition \ref{Mirkovic}.

\smallskip

\noindent (5).  We may write $\lambda = \lambda_1 - \lambda_2$, where $\lambda_i$ is antidominant ($i=1,2$), $\langle \alpha, \lambda_1 \rangle = -1$, and $\langle \alpha, \lambda_2 \rangle = 0$.  

As in the proof of (4) above, we note that $J_{s*}(J_{\lambda_1!}J_{-\lambda_2*})$ and $(J_{\lambda_1!}J_{-\lambda_2*})J_{s*}$ are each perverse, so by associativity we may unambiguously write 
$$
J_{s*}\Xi^-_{\lambda} J_{s*} = J_{s*}(J_{\lambda_1!}J_{-\lambda_2*})J_{s*}.
$$  
Using (4) and associativity, this is $(J_{s*}J_{\lambda_1!}J_{s*})J_{-\lambda_2*} = (J_{s*}\Xi^-_{\lambda_1}J_{s*})\Xi^-_{-\lambda_2}$.  (Note $J_{s*}J_{\lambda_1!}J_{s_*}$ is unambiguous and perverse, since $st_{\lambda_1} < t_{\lambda_1}$ implies that $J_{\lambda_1!}J_{s*} = J_{s!}J_{s t_{\lambda_1}!}J_{s*}$, and therefore $J_{s*}J_{\lambda_1!}J_{s*} = J_{s t_{\lambda_1}!}J_{s*}$.)

Since by (3) 
$\Xi^-_{s\lambda_1}\Xi^-_{-\lambda_2} = \Xi^-_{s\lambda}$, we are reduced to proving $J_{s*}\Xi^-_{\lambda_1}J_{s*} = \Xi^-_{s\lambda_1}$, i.e., to prove the result for general $\lambda$ it is enough to consider $\lambda$ which are antidominant.

Therefore assume $\lambda$ is antidominant, and write $l(t_\lambda) =l$.  Following Lemma 4.4 (b) of \cite{Lus4} we see
\begin{itemize}
\item $l(t_\lambda s) = l+1$ and $l(s t_\lambda) = l-1$,
\item $\lambda + s\lambda$ is antidominant,
\item $l(t_\lambda s t_\lambda) = 2l-1$ and $l(t_\lambda s t_\lambda s) = 2l-2$; in particular $l(t_{\lambda} s t_{\lambda}) = l(t_\lambda) + l(s t_\lambda)$.\end{itemize} 

Taking these relations, the previous parts of the Lemma, Corollary \ref{equivalence}, and associativity into account, we find
\begin{align*}
J_{\lambda!}\Xi^-_{s \lambda} &= \Xi^-_{\lambda} \Xi^-_{s \lambda} \\
&= \Xi^-_{\lambda + s\lambda} \\
&= J_{t_\lambda s t_\lambda s!} \\
&= J_{t_\lambda s t_\lambda!}J_{s*} \\
&= J_{\lambda!}J_{s t_\lambda!}J_{s*} \\
&= J_{\lambda!}(J_{s*} J_{\lambda!} J_{s*}).
\end{align*}
Using Corollary \ref{equivalence} again to cancel $J_{\lambda!}$ from each side, we obtain the desired equality $J_{s*} \Xi^-_{\lambda} J_{s*} = \Xi^-_{s \lambda}$.

\end{proof}

Note that property (5) is the analogue of Bernstein's relation
$$
\tilde{T}^{-1}_s \Theta^-_\lambda \tilde{T}^{-1}_s = \Theta^-_{s\lambda},
$$
which was a main ingredient in the proof of Lemma \ref{prop3.4}.  
In fact the same argument can be applied to prove the following corollary.

\begin{corollary}\label{sheaf_prop3.4}
If $\lambda$ 
is minuscule and $t_\lambda = w^\lambda w$ as in section 3, then
$$
\Xi^-_\lambda = J_{w^\lambda!}J_{w*}.
$$
Writing $w^\lambda = t_1 \cdots t_{r-p} \tau$ and $w= s_1 \cdots s_p$ as in section 3, we have
$$
\Xi^-_\lambda = J_{t_1!}\cdots J_{t_{r-p}!} J_{\tau!} J_{s_1*} \cdots J_{s_p*}.
$$

\end{corollary}

\begin{proof}
Suppose $\lambda$ is in the Weyl orbit of an antidominant minuscule coweight $\mu^-$.  We have
$$
\Xi^-_{\mu^-} = J_{\mu^-!}.
$$
Choose the sequence of simple reflections $s_1, \dots, s_p$ as in Lemma \ref{prop3.4}.  By induction on $p$, we easily see that 
$$
\Xi^-_{s_p\cdots s_1(\mu^-)} = J_{t_1\cdots t_{r-p}\tau!}J_{s_1\cdots s_p*}.
$$
Indeed, using induction and Lemma \ref{properties} this equality for $p-1$ can be written 
$$
\Xi^-_{s_{p-1} \cdots s_1(\mu^-)} = J_{s_p!}J_{t_1\cdots t_{r-p}\tau!}J_{s_1\cdots s_{p-1}*}.
$$
Multiplying on each side by $J_{s_p*}$ and using Lemma \ref{moreproperties} (5), Corollary \ref{equivalence}, and associativity yields
$$
\Xi^-_{s_p \cdots s_1(\mu^-)} = J_{t_1 \cdots t_{r-p} \tau!}J_{s_1 \cdots s_p*},
$$
as desired.  

The second statement follows from the first, using Lemma \ref{properties}.  To justify this, we need to show that
$$
(J_{t_1!}\cdots J_{t_{r-p}!}J_{\tau!})(J_{s_1*}\cdots J_{s_p*}) = J_{t_1!}\cdots J_{t_{r-p}!}J_{\tau!}J_{s_1*}\cdots J_{s_p*},
$$
where the products of the form $J_{t_1!} \cdots J_{s_p*}$ denote the $k$-fold convolution mentioned in section 6.3.  This results from the generalized associativity discussed there.

We remark that it is important here that the underlying expression $t_1 \cdots t_{r-p} \tau s_1 \cdots s_p$ is reduced.
\end{proof}

Taking Corollary 
\ref{sheaf_prop3.4} 
as well as Lemmas \ref{demazure}, \ref{properties} and \ref{moreproperties} into account, we get the following sheaf-theoretic interpretation for minimal expressions for 
$\Theta^-_\lambda$.

\begin{theorem}\label{minimal_sheaf}
\noindent {\em (a)}  Let $\lambda$ be a minuscule coweight of any root system.  Write $t_\lambda = w^\lambda w = t_1 \cdots t_{r-p} \tau s_1 \cdots s_p$, as in Lemma \ref{prop3.4}.  Then
$$
\Xi^-_\lambda = m_*({\mathcal D}),
$$
where ${\mathcal D}$ is the perverse sheaf 
$$
{\mathcal D} = (j_{t_1!}{\mathcal A}_{t_1}) \tilde{\boxtimes} \cdots 
\tilde{\boxtimes} (j_{t_{r-p}!}{\mathcal A}_{t_{r-p}}) \tilde{\boxtimes}
(j_{\tau!}{\mathcal A}_\tau) \tilde{\boxtimes} (j_{s_1*}{\mathcal A}_{s_1}) 
\tilde{\boxtimes} \cdots \tilde{\boxtimes} (j_{s_p*}{\mathcal A}_{s_p})
$$
on the Demazure resolution 
$m: X(t_1) \tilde{\times} \cdots \tilde{\times} X(t_{r-p}) \tilde{\times} 
X(\tau) 
\tilde{\times} X(s_1) \tilde{\times} \cdots \tilde{\times} X(s_p) \rightarrow X(t_\lambda)$ of the Schubert variety $X(t_\lambda)$.

\noindent {\em (b)}
Let $\lambda$ be a coweight for ${\rm GL}_n$, and write it as 
$$
\lambda = \lambda _1 + \cdots + \lambda _k,
$$
where each $\lambda_i$ is minuscule, and $l(t_\lambda) = \sum_i l(t_{\lambda_i})$.  
For each $i=1,\dots,k$ we have a decomposition and reduced expression 
$$
t_{\lambda_i} = w^{\lambda_i}w_i = t^i_1\cdots t^i_{q_i}\tau^i s^i_1 \cdots
s^i_{p_i},
$$
as in Lemma \ref{prop3.4}. 

Then $\Xi^-_{\lambda} = m_*({\mathcal D})$, where ${\mathcal D}$ is the perverse sheaf 
$$
J_{t^{1}_1!} \tilde{\boxtimes} \cdots 
\tilde{\boxtimes} J_{t^{1}_{q_1}!} \tilde{\boxtimes}
J_{\tau^{1}!} \tilde{\boxtimes} J_{s^{1}_1*} 
\tilde{\boxtimes} \cdots \tilde{\boxtimes} J_{s^{1}_{p_1}*} 
\tilde{\boxtimes} \cdots  \cdot \cdots \tilde{\boxtimes} 
J_{t^{k}_1!} \tilde{\boxtimes} \cdots  
\tilde{\boxtimes} J_{t^{k}_{q_k}!} \tilde{\boxtimes}
J_{\tau^{k}!} \tilde{\boxtimes} J_{s^{k}_1*} 
\tilde{\boxtimes} \cdots \tilde{\boxtimes} J_{s^{k}_{p_k}*} 
$$
on the Demazure resolution $m: \tilde{X}(t_\lambda) \rightarrow X(t_\lambda)$ corresponding to the reduced expression
$$
t_\lambda = (t^{1}_1 \cdots t^{1}_{q_1} \tau^{1} s^{1}_1 \cdots 
s^{1}_{p_1}) \cdots (t^{k}_1 \cdots t^{k}_{q_k} \tau^{k} s^{k}_{1} \cdots s^{k}_{p_k}).
$$

Consequently, $\Xi^-_\lambda$ is the push-forward of an explicit perverse sheaf on a Demazure resolution of $X(t_\lambda)$, for every coweight 
$\lambda$ of ${\rm GL}_n$.
\end{theorem} 

\begin{proof}
\noindent (a).  This follows directly from Corollary \ref{sheaf_prop3.4} and the definition of the $r$-fold convolution product.

\noindent (b).  By Lemma \ref{moreproperties} (3) and generalized associativity, we have
$$
\Xi^-_\lambda = \Xi^-_{\lambda_1} \star \cdots \star \Xi^-_{\lambda_k},
$$
where the right hand side denotes the $k$-fold convolution product.    

Part (b) then follows from part (a) and another application of generalized associativity.

\end{proof}

For $x \leq t_\lambda$, write $\Theta^-_\lambda(x)$ for the coefficient of $T_x$.

\begin{corollary}\label{generalformula}
Let $\lambda$ be a minuscule coweight of a root system, or an arbitrary coweight for ${\rm GL}_n$.  Let $m, {\mathcal D}$ be as in Theorem \ref{minimal_sheaf}.  
Then
$$
\Theta^-_\lambda(x) = \varepsilon_\lambda {\rm Tr}({\rm Fr}_q, H^\bullet(m^{-1}(x),{\mathcal D})),
$$
where the right hand side is the alternating trace of Frobenius on the \'{e}tale cohomology of the fiber over $x$ with coefficients in ${\mathcal D}$.
\end{corollary}

\section{Acknowledgements}

We thank I. Mirkovic for kindly allowing us to include his unpublished result 
(Proposition 6.2), which was crucial to the last section of this paper.  We are also 
grateful to R. Bezrukavnikov for some very helpful 
discussions concerning Mirkovic's theorem.  

Some of this paper was written during the first author's 2000-2001 visit to the Institute for Advanced Study in Princeton, which he thanks for support and hospitality.  His research is partially supported by an NSERC research grant and by NSF grant DMS 97-29992.

The second author was supported by an NSERC Undergraduate Student Research 
Award during the summers of 2000 and 2001. 
 
We wish to thank M. Rapoport for his comments on this paper.  We are also grateful to R. Kottwitz for his careful reading and for pointing out some errors in the first version of this paper.  We thank him for several other mathematical and expositional suggestions, in particular concerning the proof of Corollary 6.3.

\small
\bigskip
\obeylines
\noindent
University of Toronto
Department of Mathematics
100 St. George Street
Toronto, ON M5S 3G3, Canada
email: haines@math.toronto.edu

\bigskip
\obeylines
University of Chicago
Department of Mathematics
5734 S. University Ave.
Chicago, IL 60637
email: alexandra@math.uchicago.edu

\end{document}